\documentclass{amsart}

\usepackage{amssymb, amsmath, latexsym, a4}
\usepackage[latin1]{inputenc}
\usepackage[all]{xy}

\def\B{\bigskip \par}

\def \A{\mathcal A}
\def \B{\mathcal  B}
\def \P{\mathcal  P}
\newcommand{\Z}{\mathbb{Z}}

\def \Coker{\mathop{\sf Coker}\nolimits}
\def \Hom{\mathop{\sf Hom}\nolimits}

\newtheorem{De}{Definition}[section]
\newtheorem{Th}[De]{Theorem}
\newtheorem{Pro}[De]{Proposition}
\newtheorem{Le}[De]{Lemma}
\newtheorem{Co}[De]{Corollary}
\newtheorem{Rem}[De]{Remark}

\begin{document}
\title{Goodwillie's calculus via relative homological algebra. The abelian case}
\bigskip

\author[T. Pirashvili]{Teimuraz  Pirashvili}

\address{Department of Mathematics\\
University of Leicester\\
University Road\\
Leicester\\
LE1 7RH, UK} \email{tp59-at-le.ac.uk}

\maketitle

\section{Introduction} We will explain how  elementary concepts of relative homological algebra yield the Taylor tower for functors from pointed categories to abelian groups recovering the constructions of Johnson and McCarthy \cite{J2},\cite{J3}.

 Let ${\bf C, \ D}$ be abelian categories with enough projective objects. Let $i_*:{\bf C}\to {\bf D}$ and  $i^*:{\bf D}\to {\bf C}$ be functors, such that $i^*$ is left adjoint to $i_*$. We will assume that $i_*$ is full and faithful and exact. After taking the left derived functors one obtains a pair of adjoint functors
 $({\sf L}(i^*)\vdash{\sf L}(i_*))$ between the derived categories $ {\sf D}^-({\bf D})$ and ${\sf D}^-({\bf C})$. In general, ${\sf L}(i_*):{\sf D}^-({\bf C}) \to{\sf D}^-({\bf D})$ is not a full embedding. Instead one defines a full subcategory ${\sf D}^-_{\bf C}({\bf D})$ of ${\sf D}^-({\bf D})$ by
$${\sf D}^-_{{\bf C}}({\bf D})=\{X_*\in {\sf D}^-({\bf D})\,|\, H_n(X_*)\in {\bf C}, n\in \Z \}.$$
Denote by $j_*:{\sf D}^-_{\bf C}({\bf D})\to {\sf D}^-({\bf D})$ the full inclusion. Then the functor ${\sf L}(i_*)$ factors through $j_*$. In the favourable cases the  functor $j_*$ has left adjoint $j^*$, however we do not know whether $j^*$ always exists.  In the next section we will construct the functor $j^*$ under  certain circumstances. Our construction is based on the elementary results  of the relative homological algebra \cite{em} and is probably  well-known. In the last section we explain how the results  of Section \ref{2} imply the main results of \cite{J2},\cite{J3}. 

In \cite{p} we will extend our method from abelian to nonabelian case.

\section{The main construction}\label{2} Let $\A$ be an abelian category with  coproducts and let $\mathcal P$ be a set of  objects in $\A$ such that each $P\in {\mathcal P}$ is  projective. 
Define the following full subcategory 
$$\B=\P^\bot=\{A\in \A \,| \,Hom_{\A}(P,A)=0, P\in \P\}.$$ 
It is clear that $\B$ is a thick subcategory of $\A$. That is, $\B$ is closed under taking kernels, cokernels and extensions. In particular, $\B$ is also abelian. Denote by $i_*:\B\to \A$ the inclusion. Then $i_*$ is exact.

For any $A\in \A$ one puts
$$\Phi(A)=\bigoplus_{f:P\to A}P.$$
Here $P$ runs through all objects of $\P$. For a morphism $f:P\to A$ we let $in_f:P\to \Phi(A)$ be the standard inclusion. Define $\epsilon_A:\Phi(A)\to A$
by $\epsilon_A\circ in_f=f$ and  denote $\Coker(\epsilon_A)$ by $i^*(A)$. Since $\Hom_{\A}(P,\epsilon_A)$ is surjective one  sees that  $i^*(A)\in \B$. In this way one obtains a functor  $i^*:\A\to \B$ which is  left adjoint  to $i_*$.

A morphism $f:X\to Y$ in $\A$ is called $\mathcal P$-epimorphism provided $\Hom_{\A}(P,f):\Hom_{\A}(P,X)\to \Hom_{\A}(P,Y)$
is surjective. For example,  for any object $A\in \A$ the morphism $\epsilon_A:\Phi(A)\to A$  is a $\P$-epimorphism. Hence $\P$ is a projective class in the sense of \cite{em} and therefore by \cite[Proposition 3.1]{em} any object $A$ has a $\P$-projective resolution. Thus there is a chain complex $(X_*,d)$ such that $X_n=0$ if $n<-1$, $X_{-1}=A$, $X_n\in \P$ for any $n\geq 0$ and for any $P\in \P$ the following sequence is exact:
$$\cdots \to \Hom_{\A}(P,X_n)\to  \cdots \to \Hom_{\A}(P,X_0)\to \Hom_{\A}(P,X_{-1})\to 0.$$
It follows that $X_*\in {\sf D}_{\B}^-(\A).$ By the standard properties of  $\P$-projective resolutions the assignment $A\mapsto X$ extends to a functor $j^*:{\sf D}^-(\A)\to {\sf D}^-_\B(\A)$ which turns to be  left adjoint to $j_*$. 
%One  puts, $\theta_n(A)=H_n(X_*)$, $n\geq -1$. Thus  in dimensions $n\geq 1$, $\theta_n$ are the left $\P$-derived functors of the identity functor.  In practice $\theta_*(A)$ is a very interesting invariant of $A$ (and $\P$).  For examplle, if the notations %of the next section one takes  ${\bf M}$ to be the category of finitely generated free abelian groups and $n=1$, then $\theta_*$  is so called Dold-Puppe stable derived functors \cite{DP}.

Assume now that instead of a single set $\P$, a descending sequence of sets
$$ \cdots \subset \P_n\subset \P_{n-1}\subset \cdots \subset \P_1$$
is given, each of which satisfies the assumptions made in the beginning of Section \ref{2}. One obtains abelian categories $\B_n=\P_n^\bot$ and functors $i_{n*}, i^*_n,j_{n*}, j^*_n$. Clearly, $\B_1\subset \B_2\subset \B_3\subset \cdots \subset \A$ and for any object $A\in \A$ one obtains the towers of epimorphisms
$$A\to  \cdots \to i_* i^*_n(A)\to i_*i^*_{n-1}(A)\to \cdots\to  i_*i^*_2(A)\to i_*i^*_1(A)$$
and of morphisms in ${\sf D}^-(\A)$
$$ A\to \cdots \to j_*j^*_n(A)\to j_*j^*_{n-1}(A)\to \cdots \to j_*j^*_2(A)\to j_*j^*_1(A).$$

\section{Applications to Goodwillie's calculus}
Let ${\bf M}$ be  a small category with zero object $0$ and finite coproduct $\vee$. We let $\A$ be the category of all functors from ${\bf M}$ to the category of abelian groups. Then $\A$ is an abelian category with enough projective objects. The functors $h_a$ are small projective generators of $\A$. Here $a$ is running through all objects of the category ${\bf M}$ and $h_a\in \A$ is given by $h_a= \Z[\Hom_{\bf M}(a,-)]$.  The obvious maps $a\to 0\to a$ yield a splitting $h_a=\bar{h}_a\oplus \Z$, where $\Z=h_0$ is the constant functor with values equal to $\Z$. Thus the collections $\bar{h}_a$, $a\in {\bf M}$ together with $\Z$ also form a family of small projective generators.
Clearly $h_{a\vee b}=h_a\otimes h_b$. It follows that the level-wise tensor product of projective objects is again a projective object. For any natural number $n\geq 1$ we let $\P_n$ be the collection of projective objects of the form
$\bar{h}_{a_1}\otimes \cdots \otimes \bar{h}_{a_{k}}$, $k>n$. One easily checks that the corresponding category $B_n=\P_n^\bot$ is the category of functors  of  degree $\leq n$ (in the sense of Eilenberg-MacLane), while ${\sf D}^-_{\B_n}(\A)$ is equivalent to the category of functors from $\bf M$ to the category of chain complexes of abelian groups of degree $\leq n$ (in the sense of Goodwillie). This follows from the fact that $\Hom_\A(\bar{h}_{a_1}\otimes\cdots \otimes \bar{ h}_{a_k}), T)=cr_kT(a_1\cdots,a_k)$, where $cr_k$ is the $k$-th crossed-effect \cite{J3}. The last isomorphism is a trivial consequence of the Yoneda lemma and the decomposition rule: $h_{a\vee b}=h_a\otimes h_b$.
It follows that in this situation the towers constructed in Section \ref{2} and the ones constructed in \cite{J2},\cite{J3} are equivalent.

\end{document}